\documentclass{amsart}
\usepackage{amscd,amsmath,amssymb,amsfonts,verbatim, xcolor,graphicx}
\usepackage{calrsfs}
\usepackage{latexsym, xcolor, amsmath, amscd}
\usepackage[all]{xy}
\DeclareMathAlphabet{\pazocal}{OMS}{zplm}{m}{n}
\allowdisplaybreaks
\usepackage{calrsfs}
\makeatletter
\newcommand*\dotp{\mathpalette\dotp@{.5}}
\newcommand*\dotp@[2]{\mathbin{\vcenter{\hbox{\scalebox{#2}{$\m@th#1\bullet$}}}}}
\makeatother

\makeatletter
\newcommand*\bigcdot{\mathpalette\bigcdot@{.5}}
\newcommand*\bigcdot@[2]{\mathbin{\vcenter{\hbox{\scalebox{#2}{$\m@th#1\bullet$}}}}}
\makeatother

\theoremstyle{definition}

\newtheorem*{theorem*}{Theorem}
\newtheorem*{conjecture*}{Conjecture}
\newtheorem*{remark*}{Remark}

\numberwithin{equation}{subsection}
\title{A survey of the additive dilogarithm}

 \subjclass[2010]{19E15, 14C25}
\begin{document}
\author{S\.{I}nan \"{U}nver}
\address{Ko\c{c} University, Mathematics Department. Rumelifeneri Yolu, 34450, Istanbul, Turkey}
\email{sunver@ku.edu.tr}
\maketitle
\noindent

\begin{abstract}
 Borel's construction of the regulator gives an injective map from the algebraic $K$-groups of a number field to its Deligne-Beilinson cohomology groups. This has  many interesting arithmetic and geometric consequences. The formula for the regulator is expressed in terms of the classical polyogarithm functions. In this paper, we give a survey of the additive dilogarithm and the several different versions of the weight two regulator in the  infinitesimal setting. We follow a historical approach which we hope will provide motivation for the  definitions and the constructions. 
\end{abstract}

\section{Introduction} 

The dilogarithm function, even though it has been known for a very long time, has  become more prevalent in the past few decades because of its relation to  regulators in algebraic $K$-theory, as was first observed in the pioneering work of Bloch \cite{bloch-reg}.  Among others, this point of view was furthered through the far reaching conjectures of  Beilinson on motivic cohomology  \cite{beil}, by the work of Zagier on his conjecture relating special values of Dedekind zeta functions of number fields to values of regulators \cite{zag-pol} and in many works of Goncharov (\cite{gon2}, \cite{vol}, \cite{gon3} to name a few). The dilogarithm function also appears in hyperbolic geometry, conformal field theory and the theory of cluster algebras. The  survey \cite{zagier} is an excellent introduction to some aspects of this function. 

In this note, we give a survey of the infinitesimal version of the above theory. Since the generalizations of the results in this survey to  higher weights is still in progress, we restrict to the case of  the dilogarithm. In \textsection \ref{higher weight}, we will only briefly mention the construction of    additive polylogarithms of higher weight on certain special linear configurations. The existence of this theory itself is quite surprising and is based on ideas of Cathelineau (\cite{cat-hom}, \cite{cat}), Bloch and Esnault \cite{BE2} and Goncharov \cite{vol}, which we will describe in detail below. We emphasize that these functions  cannot be deduced from their classical counterparts through a limiting process. We  illustrate this point in the, somewhat deceptively simple, case of weight 1 as follows. The regulator over the complex numbers is given essentially by  the real analytic map $\log |\cdot|: \mathbb{C} ^{\times} \to \mathbb{R}.$ On the other hand, in the infinitesimal case, for $k$ a field of characteristic 0, and $k_{n}:=k[[t]]/(t^n),$ one  has the algebraic map  $\log ^{\circ}: k_n ^{\times } \to k_n$ defined by $\log ^{\circ} (a):=\log (\frac{a}{a(0)}).$  The use of the absolute value makes the first function non-algebraic, single valued  and dependent, in an essential way, on the  local field in question.  In the second case, the map $k \to k_n,$ which is a section of the canonical projection from $k_n$ to its quotient by its nilradical achieves the  purpose of choosing a branch in an appropriate sense. We will see below that over a scheme with non-reduced structure such local splittings, which correspond to retractions of the scheme with the reduced induced structure, will play a role analogous to choosing  branches.

In the second section, we briefly recall the definitions of the Bloch-Wigner dilogarithm, the Chow dilogarithm of Goncharov and  Bloch's regulator function  from $K_{2}(C)^{(2)} _{\mathbb{Q}}$ of a curve $C.$ We emphasize the point of view of the Aomoto dilogarithms and scissors congruence class groups whose  analogs will be the main motivation for the infinitesimal versions of the above functions.

In the third section, we give the infinitesimal analogs of these functions. Starting with the ideas of Cathelineau, Goncharov and Bloch-Esnault. We also recall the additive dilogarithm construction of Bloch-Esnault.  

In the fourth section, we discuss the construction of the infinitesimal Chow dilogarithm, together with its application to algebraic cycles and Goncharov's strong reciprocity conjecture. We also describe  the infinitesimal version of Bloch's regulator on curves. 

In the last section, we discuss some partial results in higher weights and in characteristic $p$ and some open problems.

{\bf Conventions.} Except in \textsection \ref{char p}, we will consider motivic cohomology always with $\mathbb{Q}$-coefficients. Therefore all the Bloch complexes, Aomoto complexes etc. are tensored with $\mathbb{Q}.$ For example, the notation $\Lambda ^2 k^{\times}$ means that the  group $\Lambda ^{2} _{\mathbb{Z}} k^{\times}$ is tensored with $\mathbb{Q}.$ The cyclic homology and Andr\'{e}-Quillen homology groups are always considered relative to $\mathbb{Q}.$ The notation $\Omega^1 _{A}$ for an algebra $A$ over a field always means the K\"{a}hler differentials relative to the prime field.  For an $A$-module $I,$  $S^{\bigcdot} _{A}  I$ denotes the symmetric algebra of $M$ over $A.$  For a ring $A,$ $A^{\flat}$ denotes the set of all units $a$ in $A$ such that $1-a$ is also a unit.  For a functor $F$ from the category of  pairs $(R,I)$ of rings $R$ and nilpotent ideals $I$ to an abelian category,  we let $F^{\circ}(R,I)$ denote the kernel of the map from $F(R,I)$ to $F(R/I,0).$ We informally refer to this object as the {\it infinitesimal part of} $F.$ We have the corresponding notion for the category of artin local algebras over a field, since their maximal ideals are nilpotent.   

\section{Bloch-Wigner dilogarithm and the scissors congruence class group}

\subsection{Aomoto dilogarithm}\label{aomoto}
The general conjectures on motives expect that for any field $k$ one has a tannakian category ${\rm MTM}_{k}$ over $\mathbb{Q}$   of mixed Tate motives over $k.$ This gives a  graded Hopf algebra $\mathcal{A}_{\bigcdot}(k)$ such that a mixed Tate motive over $k$  is the same as a graded  $\mathbb{Q}$-space with a co-module structure over $\mathcal{A}_{\bigcdot}(k).$

 Since the objects in ${\rm MTM}_{k}$ should be constructed  from Tate objects by means of extensions, one expects   ${\rm MTM}_{k}$ to  have a {\it linear} algebraic description. In \cite{bgsv} a  graded Hopf algebra $A_{\bigcdot}(k)$ was defined, using linear algebraic objects, such that one expects a natural map  $A_{\bigcdot}(k) \to \mathcal{A}_{\bigcdot}(k).$  

This $A_{\bigcdot} (k)$ is the graded Hopf algebra of Aomoto polylogarithms over $k$ defined in \cite[\textsection 2]{bgsv}. An $n$-simplex $L$ in $\mathbb{P}^n _{k}$ is an $(n+1)$-tuple $(L_{0},\cdots,L_{n})$ of hyperplanes. It is said to be non-degenerate if the hyperplanes are in general position. A pair of simplices $(L,M)$ is said to be admissable if they do not have a common face.  $A_{n}(k)$ is the $\mathbb{Q}$-space generated  by pairs of admissable simplices $(L;M)$ in $\mathbb{P}^{n} _{k}$ subject to the following  relations: 

(i) $(L,M)=0,$ if one of the simplices is degenerate; 

(ii) $(L,M)$ is anti-symmetric with respect to the ordering of the hyperplanes in both of the $n$-simplices. 

(iii) If $L$ is an $n+2$-tuple of hyperplanes $(L_0,\cdots, L_{n+1})$ and $\hat{L}^j$ is the $n$-simplex obtained by omitting $L_{j},$ then 
$
\sum _{0 \leq j \leq n+1} (\hat{L}^j,M)=0,
$ and the corresponding relation for the second component. 

(iv) For  $\alpha \in {\rm GL}_{n+1}(k),$ $(\alpha(L),\alpha(M))=(L,M).$

There are certain configurations, called polylogarithmic configurations, in $A_{n}(k)$ that play an important role in understanding the motivic cohomology of $k,$ since they act as building blocks for all configurations \cite[\textsection 1.16]{bgsv}. Let $P_{n}(k)$ denote the subgroup of {\it prisms} in $A_{n}(k).$ This is the subgroup generated by configurations which come from products of configurations from lower dimensions. For every $a \in k^{\flat}:=k^{\times}\setminus \{1 \},$ there is a special   configuration $(L,M_{a}) \in A_{n}(k)$ \cite[Fig. 1.14]{gon2}, which corresponds to the value of the abstract polylogarithm at $a.$ If $z_{i},$ $0\leq i \leq n$ are the homogenous coordinates on $\mathbb{P}^n _{k},$ then $L_{i}$ is defined by $z_{i}=0.$ The simplex $M_a$ is defined by the following formulas. $M_0: z_{0}=z_1;$ $M_{1}: z_{0}=z_1+z_2;$ $M_{i}:z_i=z_{i+1}, $ for $2\leq i <n;$ and $M_n:az_0=z_n.$

This defines a map $
l_{n}:\mathbb{Q}[k^{\flat}] \to A_{n}(k)/P_{n}(k),
$ which sends the generator $[a]$ to the class of $(L,M_a).$
 Denoting the image of $l_{n}$  by $B_{n} '(k),$ one expects the co-multiplication on $A_{\bigcdot}(k)$ to induce a complex $\Gamma_{k}'(n)$:
$$
B_{n} '(k) \to B_{n-1} '(k) \otimes k^{\times}\to  \cdots \to B_{2} '(k) \otimes \Lambda^{n-2}k^{\times} \to \Lambda^{n} k^{\times}, 
$$
which would compute the motivic cohomology of $k$ of weight $n.$

For $n=2,$ there is a simpler complex, namely the Bloch complex $\Gamma_{k}(2)$ of weight two, which computes the motivic cohomology. Let $B_{2}(k)$    be the quotient of $\mathbb{Q}[k^{\flat}],$ the vector space with basis $[x]$ for   $x \in k^{\flat},$ by the subspace generated by elements  of the form 
 \begin{eqnarray}\label{funct-dilog}
 [x]-[y]+[y/x]-[(1-x^{-1})/(1-y^{-1})]+[(1-x)/(1-y)],
 \end{eqnarray}
 for all $x,y \in k^{\times} $ such that $(1-x)(1-y)(1-x/y) \in k^{\times}.$ The last equation is the 5-term functional equation of the dilogarithm. 
Let  $\delta$ be  the map that sends $[x]$ to $(1-x)\wedge x \in \Lambda ^{2}k^{\times}.$ This map factors through $B_{2}(k)$ and we obtain a complex:
$$
\xymatrix{
B_{2}(k) \ar^{\delta}[r]& \Lambda^2 k^{\times},}
$$
concentrated in degrees 1 and 2. We denote this complex by $\Gamma_{k}(2).$ 
This complex indeed computes the motivic cohomology of $k$ with coefficients $\mathbb{Q}(2),$ by a theorem of Bloch. In other words, the sequence 
$$
0\to K_{3}(k)^{(2)} _{\mathbb{Q}} \to B_{2}(k) \to \Lambda^{2} k^{\times} \to K_{2} ^{M}(k)_{\mathbb{Q}}
 \to 0$$
 is exact (\cite{bloch-reg}, \cite{sus}).

The map $l_2$ factors through the quotient $\mathbb{Q}[k^{\flat}] \to B_{2}(k)$ to induce an isomorphism:
$$
l_2: B_{2}(k) \to B_{2 } '(k)=A_{2}(k)/P_{2}(k)
$$
which we continue to denote with the same symbol 
\cite[Proposition 3.7]{bgsv}.
This can be thought of as the abstract motivic dilogarithm function.

\subsection{Bloch-Wigner dilogarithm}
The $n$-th polylogarithm function is defined inductively by $\ell i_{1}(z)=-\log (1-z)$ and
$$
d \ell i_{k}(z)=\ell i_{k-1}(z)\frac{dz}{z}, 
$$
with  $\ell i_{k}(0)=0.$ These functions have the power series expansion $\ell i_{k}(z)=\sum _{1 \leq n}\frac{z^n}{n^k},$ in the unit disc around 0, and have multi-valued analytic continuations to $\mathbb{C} ^{\times} \setminus \{1\}.$ They appear as coordinates of a matrix which describe a canonical quotient of the fundamental groupoid associated to the Hodge realization of the  unipotent fundamental group of $\mathbb{P}^{1} \setminus \{0, 1, \infty \}$ \cite{be-de}. The specialization of this construction at a point $a \in \mathbb{C}^{\flat}$ gives a motive which coincides with the motive associated to the configuration $l_2(a)$ in \textsection \ref{aomoto}. 

The Hodge realization of this motive (specialized at a point) as well as of the motive above defined by the configurations in \textsection \ref{aomoto} above are Hodge-Tate structures. An $\mathbb{R}$-Hodge-Tate structure is a mixed $\mathbb{R}$-Hodge structure such that for every $r \in \mathbb{Z},$ its graded piece of degree $-2r$ with respect to the weight filtration are direct sums of the Tate  structures $\mathbb{R}(r),$ of weight $-2r;$  and its graded pieces of odd degree are equal to 0.   Let $\mathcal{H}_{ \bigcdot}$ denote the graded Hopf algebra associated to the tannakian category of  $\mathbb{R}$-Hodge-Tate structures.  The Hodge realization functor should give a  morphism $\mathcal{A}_{\bigcdot}(\mathbb{C}) \to \mathcal{H}_{ \bigcdot}$ of graded Hopf algebras.   

 A construction of Beilinson and Deligne (\textsection 2.5, \cite{be-de}; pp. 248-249, \cite{gon2}) associates to each framed $\mathbb{R}$-Hodge-Tate structure a number.  Associated to the variation of Hodge structures on $\mathbb{G}_{m}$ that gives the function $\log(z)$ one gets the corresponding single valued function $\log |z|.$  This  construction gives a map $p_{\mathcal{H},n}:\mathcal{H}_{n} \to \mathbb{R}.$ It turns out that this map vanishes on the products \cite{gon2}. Hence composing with the Hodge realization map associated to the Aomoto configurations, the corresponding map vanishes on prisms and one gets a map $B_{2}'(\mathbb{C})=A_{2}(\mathbb{C})/P_{2}(\mathbb{C}) \to \mathbb{R}.$ The composition of this map with $l_2: B_{2}(\mathbb{C}) \to B_{2}'(\mathbb{C})$ turns out to be, up to scaling, the Bloch-Wigner dilogarithm $D$ defined by 
 $$
D(z)={\rm Im}(\ell i_{2}(z))+{\rm arg}(1-z) \log |z|.
 $$
 The main importance of the Bloch-Wigner dilogarithm comes from the fact that they are regulators. 
 
 Composing $p_{\mathcal{H},n}$  with the Hodge realization would give a map   
 $$
{\rm vol}_{n}: \mathcal{A}_{n}(\mathbb{C}) \to \mathbb{R},
 $$
which is an analog  of the volume map on the scissors congruence class groups below and its infinitesimal version is the main concern of this survey.

\subsection{Chow dilogarithm}\label{section chow dilog}

If $X/\mathbb{C}$ is a smooth and projective curve over $\mathbb{C},$ there is a version of the dilogarithm above which gives  certain regulators of $X.$ Namely ${\rm H}^{3} _{\pazocal{M}}(X,\mathbb{Q}(3))\simeq K_{3}(X)^{(3)} _{\mathbb{Q}} $ and applying the Leray-Serre spectral sequence to the map $X \to \mathbb{C},$ there would be a map 
$K_{3}(X)^{(3)} _{\mathbb{Q}}\simeq {\rm H}^{3} _{\pazocal{M}}(X,\mathbb{Q}(3)) \to {\rm H}^{1} _{\pazocal{M}}(\mathbb{C},H^{2}(X/\mathbb{C})(3) )={\rm H}^{1} _{\pazocal{M}}(\mathbb{C},\mathbb{Q}(2))\simeq K_{3}(\mathbb{C})_{\mathbb{Q}} ^{(2)}.$ Combining with the regulator $K_{3}(\mathbb{C})_{\mathbb{Q}} ^{(2)}\to \mathbb{R},$ given by the Bloch-Wigner dilogarithm above, one would get a map 
$$
K_{3}(X)^{(3)} \to \mathbb{R}.
$$
This map is  given by the following Chow dilogarithm of Goncharov.  

If  $f_{1},\, f_2, \,$ and $f_3$  are rational functions on $X.$ Let
$$
r_{2}(f_{1},f_{2},f_{3}):={\rm Alt}_{3}(\frac{1}{6} \log |f_1| \cdot d \log |f_{2}| \wedge d \log |f_{3}| -\frac{1}{2}\log |f_1|\cdot d \arg f_{2} \wedge d \arg f_{3} ),
$$
which has the formal property that $d(r_{2}(f_{1},f_{2},f_{3}))={\rm Re} (dlog (f_1)\wedge dlog (f_3)\wedge dlog (f_3)).$ The map 
$
\rho_{\mathbb{R}}: \Lambda^{3} \mathbb{C}(X) ^{\times} \to \mathbb{R}
,$ given by 
\begin{eqnarray*}\label{realchow}
\rho_{\mathbb{R}}(f_1 \wedge f_2 \wedge f_3) := \int_{X(\mathbb{C})} r_{2}(f_{1},f_{2},f_{3}),
\end{eqnarray*}
is, up to a constant multiple, the Chow dilogarithm  \cite[p. 4]{gon3}. The middle cohomology of  the  complex  
$$
\cdots \to B_{2}(\mathbb{C}(X)) \to (\oplus _{x \in X}B_{2}(\mathbb{C})) \oplus \Lambda^{3} \mathbb{C}(X) ^{\times} \to \oplus _{x \in X} \Lambda ^2\mathbb{C}^{\times}\to \cdots 
$$
is $K_{3}(X)^{(3)} _{\mathbb{Q}}$ and the map $(\oplus _{x \in X}D ) \oplus \rho_{\mathbb{R}}$ obtained  by using the Bloch-Wigner and the Chow dilogarithm, gives the regulator.

\subsection{Bloch's regulator on curves and the tame symbol}\label{section bloch regulator} 

There is another regulator which is based on a version of the dilogarithm. Again assume that $X/\mathbb{C}$ is a smooth and projective curve. This regulator is essentially the map from $K_{2}(X)^{(2)}_{\mathbb{Q}}$ to the corresponding  Deligne cohomology group:  $K_{2}(X)^{(2)} _{\mathbb{Q}} \to {\rm H}^{2} _{D}(X_{an},\mathbb{Q}(2))\simeq H^{1} (X_{an},\mathbb{C}/\mathbb{Q}(2)).$ Dividing by $2 \pi i$ and using the exponential map on the coefficients, the last cohomology group is identified with $H^{1}(X_{an},\mathbb{C}^{\times} _{\mathbb{Q}}).$  Since $H^{1}(X_{an},\mathbb{C}^{\times})$ coincides with local systems of complex vector spaces of rank 1 and hence with analytic line bundles with connection.   The above  map can also be deduced from the local and analytic  construction of Deligne, which associates to each pair $f,g$ of meromorphic functions on $X,$ a line bundle with connection on $X_{an},$  such that the monodromy at each point is given by the tame symbol of $f$ and $g$ at that point \cite{de}. Explicitly, if $\log (f)$ is a  choice of a branch of $f,$ locally analytically, then the line bundle in question is the trivial line bundle with the connection $\nabla$ given by $\nabla(1)=\frac{1}{2\pi i} \log (f)\frac{dg}{g}.$ For a different choice $\log(f)+n2 \pi i$ of a logarithm of $f,$ the isomorphism between the line bundles with connection is given as multiplication by $g^{-n}$ \cite[\textsection 2.3]{de}.

When $X$ is defined over a number field, the Bloch regulator is fundamental in the study of certain special values of the $L$-function of $X$ \cite{ram}. It also appears, for example, in the geometric study of cycles on $X/\mathbb{C}$ \cite{gg}.

\section{Additive dilogarithm and the infinitesimal scissors congruence class group}  

In this section, we start with the 4-term functional equation for the entropy function which is also satisfied by an infinitesimal version of the Dehn invariant for scissors congruence class groups. This 4-term functional equation of Cathelineau can be thought of as a deformation of the 5-term functional equation that is restricted to certain special elements. The precise relation is explained in \textsection \ref{section comparison}. Next we describe Goncharov's idea that the hyperbolic scissors congruence class group can be thought of degenerating to the euclidean one as the model for the hyperbolic space blows up. We continue the section with describing the construction of the additive dilogarithm by Bloch and Esnault based on the localization sequence in $K$-theory and end the section on  our construction of the additive dilogarithm on the Bloch group.

\subsection{The 4-term functional equation}\label{4-term-section}

In information theory, Shannon's binary entropy function $H$ is defined as  
$$
H(p):=-p
\log (p)-(1-p)  \log (1-p),
$$ 
for the probability $p.$  This function satisfies the following fundamental functional equation of information theory:  
\begin{align}\label{shannon}
  H(p) + (1 - p)H (\frac{q}{1-p}) = H(q) + (1 - q)H (\frac{p}{1-q}).  
\end{align}

The same functional equation reappeared in \cite{cat-hom} as follows. For a field $k$ of characteristic 0, let $\beta_{2}(k)$ is the vector space over $k$ generated by the symbols $\langle a \rangle ,$ for $a \in k^{\flat}$ with relations generated by 
\begin{align}\label{add-func-eq}
    \langle p \rangle -\langle q \rangle+p\langle \frac{q}{p}\rangle + (1-p)\langle  \frac{1-q}{1-p} \rangle=0
\end{align}
when $p \neq q.$ 
These relations already imply that $\langle p \rangle= \langle 1-p \rangle $ and $\langle \frac{1}{p} \rangle=-\frac{1}{p}\langle p \rangle$ and using these the two relations (\ref{shannon}) and (\ref{add-func-eq}) are equivalent. In \cite[Theor\`eme 1]{cat-hom}, Cathelineau proves that the following sequence 
$$
\xymatrix{
0 \ar[r] & \beta_{2}(k) \ar^{D}[r]&  k \otimes k^{\times} \ar^{L}[r]&  \Omega ^{1} _{k}  \ar[r] & 0
}
$$
is exact, where $D$ is defined on the generators by $D(\langle a \rangle):=a \otimes a +(1-a) \otimes (1-a)$ and $L$ sends $a \otimes b$ to $a  \frac{db}{b}.$ 
This was used in \cite{cat-hom} in order to show that for an algebraically closed field $k$ of characteristic 0, the homology groups of ${\rm SL}(2,k)$ with adjoint action on its Lie algebra $\mathfrak{sl}(2,k)$ are given by: 
\begin{align*}
   &  {\rm H}_{1}({\rm SL}(2,k),\mathfrak{sl}(2,k))\simeq \Omega^1 _{k}\\
   & {\rm H}_{2}({\rm SL}(2,k),\mathfrak{sl}(2,k))=0
\end{align*}
This in analogy with the computation of the homology of the discrete special orthogonal group ${\rm SO}^{\delta}(3,\mathbb{R})$  with the standard 
action on $\mathbb{R}^3:$ 
\begin{align*}
   &  {\rm H}_{1}({\rm SO}^{\delta}(3,\mathbb{R}),\mathbb{R}^3)\simeq \Omega^1 _{\mathbb{R}}\\
   & {\rm H}_{2}({\rm SO}^{\delta}(3,\mathbb{R}),\mathbb{R}^3)=0.
\end{align*}
This is a restatement of Sydler's theorem that the Dehn invariant and the volume completely determine the scissors congruence class. In this euclidean case, the analog of $k\otimes k^{\times}$ is the group $\mathbb{R}\otimes \mathbb{R}/\pi \mathbb{Z}$ and the analog of the map $L$ above is   the map 
$$
\mathbb{R} \otimes \mathbb{R}/\pi \mathbb{Z} \to \Omega^1 _{\mathbb{R}} 
$$
that sends $l \otimes \theta$ to $l \frac{d(cos \,\theta )}{sin\, \theta}. $

The above $D$ can be thought of as the infinitesimal version of the Dehn invariant and the functional equation above can  be thought of as the infinitesimal version of the functional equation of the dilogarithm in the following sense.

\subsection{Hyperbolic space degenerating to  euclidean space}

In this section, we describe how  Goncharov's idea on the degeneration of  hyperbolic space to euclidean space and the analogy between the scissors congruence class groups and mixed Tate motives leads one to expect a volume map on mixed Tate motives over dual numbers which is reminiscent of the polylogarithm functions. 

\subsubsection{} \label{section scissors} If $\mathcal{G}^{n}$ is one of the three $n$-dimensional classical geometries: $\mathcal{E}^{n},$ the euclidean; $\mathcal{H}^{n},$ the hyperbolic; or $\mathcal{S}^{n},$ the spherical, then let $\mathcal{P}(\mathcal{G}^{n})$ denote the scissors congruence class group corresponding to $\mathcal{G}^{n}$. The Dehn invariant map :
$$
D_{n} ^{\mathcal{G}}: \mathcal{P}(\mathcal{G}^{n}) \to \oplus_{i=1} ^{n-2} \mathcal{P}(\mathcal{G}^{i}) \otimes \mathcal{P}(\mathcal{S}^{n-i-1})
$$
endows $\oplus \mathcal{P}(\mathcal{S} ^{\bigcdot})$ with the structure of a co-algebra and, $\oplus \mathcal{P}(\mathcal{H}^{\cdot})$ and $\oplus \mathcal{P}(\mathcal{E}^{\cdot})$  with structures of co-modules over this co-algebra \cite{vol}. 

There exists a   map from  $\mathcal{P}(\mathcal{H}^{2n-1})$ to $\mathcal{A}_{n}(\mathbb{C}),$ defined by Goncharov, which attaches a framed mixed Tate motive to an element in the  hyperbolic scissors congruence class group \cite{vol}. If one considers the Cayley spherical model for the hyperbolic geometry then as the sphere gets bigger the hyperbolic geometry approaches  the euclidean geometry \cite{vol}. Therefore, in the limit case one would expect to have a  map 
$
\mathcal{P}(\mathcal{E}^{2n-1}) \to \mathcal{A}_{n} ^{\circ}(\mathbb{C}_{2}).
$

These suggest a close similarity between  the structures of $\mathcal{A}_{n} ^\circ(k_{2})$ and $\mathcal{P}(\mathcal{E}^{2n-1} _{k})$ \cite{vol}, \cite{euc}.  The euclidean scissors congruence class group has a volume map  
$$
\mathcal{P}(\mathcal{E}^{2n-1} _{k})\to k,
$$ 
which is conjectured to induce an isomorphism from ${\rm H}^{1}(\oplus_{2n-1} \mathcal{P}(\mathcal{E} ^{\bigcdot} _{k} ) ),$ the kernel in  $\mathcal{P}(\mathcal{E}^{2n-1}_{k})$ of the Dehn invariant map,  to $k.$ For $n=2$ and $k=\mathbb{R},$ this is Sydler's theorem. In analogy, we expect a  map:
$$
{\rm vol}_{n} ^{\circ}: \mathcal{A}_{n} ^{\circ}(k_{2}) \to k,
$$
which induces an isomorphism from ${\rm H}^{1}(\mathcal{A}_{\bigcdot} ^{\circ} (k_{2})(n))$ to $k.$ Moreover, we should have the identity ${\rm vol}_{n} ^{\circ}\circ \rho_{\lambda}=\lambda^{2n-1}{\rm vol}_{n} ^{\circ},$ for $\lambda \in k ^{\times}.$ This map would be  an analog of both the map $\mathcal{A}_{n}(\mathbb{C}) \to \mathbb{R}$ that is constructed using the Beilinson-Deligne construction and of the volume map on euclidean scissors congruence class groups. 

\subsubsection{} \label{section single valued l_n} Given an element $(L,M) $ in $A_{n}(\mathbb{C}),$ this defines a framed mixed Tate motive in $\mathcal{A}_{n}(\mathbb{C})$ whose  associated mixed Hodge structure ${\rm H}^n(\mathbb{P}^n _{\mathbb{C}}\setminus L,M\setminus L ),$  is Hodge-Tate. Therefore, using the construction of Beilinson and Deligne described above, which attaches a real number to $\mathbb{R}$-Hodge-Tate structures, we get ${\rm vol}_n(L,M) \in \mathbb{R}.$  This vanishes on the products \cite{gon2} to give:
$
{\rm vol}_n: A_{n}(\mathbb{C})/P_{n}(\mathbb{C}) \to \mathbb{R}.
$
Composing with the abstract polylogarithm map induces 
$$
{\rm vol}_{n} \circ l_{n}:B_{n}'(\mathbb{C}) \to \mathbb{R}.
$$
This  has the following description. Let  $\mathcal{L}_{n}$ be the  real single valued version of the $n$-polylogarithm:
 $$
 \mathcal{L}_{n}(z):=\mathcal{R}_{n}(\sum_{j=0} ^{n} \frac{2^{j}B_{j}}{j!} (\log|z|)^{j}  \ell i_{n-j}(z)),
 $$
 where $B_{n}$ is the $n$-th Bernoulli number; $\mathcal{R}_{n}$ is the real part if $n$ is odd and the imaginary part if $n$ is even; and $\ell i _{0}(z):=-1/2.$ Then for $z \in \mathbb{C}^{\flat},$ ${\rm vol}_{n}\circ l_{n} (z)
 =\mathcal{L}_{n}(z)$ \cite{gon2}.

\subsubsection{ } Let $k$ be any field of characteristic 0. The definitions of $A_{n}(k), P_{n}(k), l_{n}$ and $B_{n} ' (k)$   exactly carry over to the  case of $k_2$ to define the groups $A_{n}(k_{2}), P_{n}(k_{2}),$ and  $B_{n} '(k_{2}),$ and a map, 
$l_{n}: \mathbb{Q}[k_{2} ^{\flat}] \to A_{n}(k_{2}).$  One would like to define a map $$
{\rm vol}_{n} ^{\circ}:A_{n}(k_{2})/P_{n}(k_2) \to k,
$$
which would be an analog of the map defined above over the complex numbers using the Beilinson-Deligne construction. This map would be the composition of the natural map from  $A_{n}(k_2).$ In this context the analog of the single valued polylogarithm $\mathcal{L}_{n}$ would be the composition ${\rm vol}_{n} ^{\circ}\circ l_{n}.$

 \subsubsection{} As in \textsection \ref{aomoto}, one has a complex  $\Gamma _{k_{2}} '(n),$ concentrated in  degrees $[1,n]:$ 
$$
B_{n} '(k_{2}) \to B_{n-1} '(k_{2}) \otimes k_{2} ^{\times}\to  \cdots \to B_{2} '(k_{2}) \otimes \Lambda^{n-2}k_{2}  ^{\times} \to \Lambda^{n} k_{2} ^{\times} 
$$
induced by the co-multiplication map on $A_{\bigcdot}(k_{2})$   and such that $\{ x\}_{i} \otimes y \in B_{i} '(k_{2}) \otimes \Lambda ^{n-i} k_{2} ^{\times}$ is mapped to: 
\begin{eqnarray*}\label{blform1}
\{ x\}_{i-1} \otimes x\wedge y \in B_{i-1} '(k_{2}) \otimes \Lambda ^{n-i+1} k_{2} ^{\times}
\end{eqnarray*}
if $i \geq 3,$ and to 
\begin{eqnarray*}\label{blform2}
(1-x)\wedge x \wedge y \in \Lambda ^{n} k_{2} ^{\times} 
\end{eqnarray*}
if $i=2.$  One expects the cohomology groups to be given by  ${\rm H}^{i}(\Gamma_{k_{2}} ' (n) )\simeq K_{2n-i}(k_{2})_{\mathbb{Q}} ^{(n)}.$ By  Goodwillie's theorem \cite{good}, we have, 
$
K_{2n-i} ^{\circ}(k_{2})_{\mathbb{Q}} ^{(n)}\simeq {\rm HC}_{2n-i-1} ^{\circ}(k_{2})^{(n-1)}.
$   The infinitesimal part of the  cyclic homology of $k_{2}$ is computed as 
$
{\rm HC}_{n} ^{\circ}(k_{2})^{(m)} \simeq \Omega ^{2m-n}  _{k},
$
for $[\frac{n+1}{2}]\leq m \leq n,$ and is 0 otherwise \cite{lambda}. Moreover, for $\lambda \in k ^{\times},$ the automorphism $\rho_{\lambda}$ of $k_{2}$ that sends $t$ to $\lambda t$ induces
multiplication by $\lambda ^{2(n-m)+1}$ on $\Omega^{2m-n} _{k}$ \cite{lambda}.   
 Combining  these, one expects the infinitesimal part of the  cohomology of $\Gamma' _{k_2}(n)$ to be: 
\begin{eqnarray*}\label{infinmotcoh}
{\rm H}^{i}(\Gamma_{k_{2}}'^{\circ}  (n)  )\simeq \Omega^{i-1} _{k}, 
\end{eqnarray*}
for $1\leq i \leq n,$ and that $\rho_{\lambda}$ induces multiplication by $\lambda^{2(n-i)+1}$ on  $\Omega^{i-1} _{k}.$ Note that when $i=1,$ this map scales  by $\lambda ^{2n-1},$ exactly like the volume map in \textsection \ref{section scissors}.

\subsection{Bloch and Esnault's construction of the additive dilogarithm on the localization sequence} The work of Bloch and Esnault was the principal motivation for the  various generalizations of the additive dilogarithm. Here we  briefly describe their work,  generalized to the case of higher moduli. The proofs of the statements can be found in \cite{BE2} and in   \cite[\textsection 6.2]{unv1}. In this section, we assume that $k$ is algebraically closed in addition to being of characteristic 0. 

Let $\pazocal{O}$ be the local ring of $\mathbb{A}^{1} _{k}$ at 0.  The localization sequence of the pair $(k[t],(t^m))$ gives the following two exact sequences: \[
\begin{CD}
K_{2}(k[t] , (t^{m})) \to  K_{2}(\pazocal{O},(t^{m})) @>{\partial}>> \oplus_{x \in k^{\times}} K_{1}(k) \to K_{1} (k[t], (t^m))\to 0
\end{CD}
\]
and 
\[
\begin{CD}
0 \to K_{1}(\pazocal{O},(t^m)) @>{\partial}>>  \oplus _{x \in k^{\times}}K_{0}(k) \to K_{0}(k[t],(t^m))\to 0.
 \end{CD}
\]
The group $$T_{m}B_{2}(k):=(K_{2}(\pazocal{O},(t^m))/im(K_{1}(k) \cdot K_{1}(\pazocal{O},(t^m)))_{\mathbb{Q}}  ,$$
is the infinitesimal analog of the Bloch group. Since $K_{0}(k[t],(t^m)) \simeq 1+(t)=(k_{m} ^{\times})^{\circ} \subseteq k_{m} ^{\times},$ the quotient   $\oplus_{x \in k^{\times}} K_{1}(k)/\partial (K_{1}(k) \cdot K_{1}(\pazocal{O},(t^m))) \simeq k^{\times} \otimes (k_{m} ^{\times})^{\circ}.$ This gives the complex : 
$$
T_{m}B_{2}(k) \to k^{\times} \otimes (k_{m} ^{\times})^{\circ},
$$
which is the analog of the Bloch complex and is denoted by $T_{m}\mathbb{Q}(2)(k).$ The cohomology groups of this complex in degrees 1 and 2 are respectively, $K_{3} ^{\circ} (k_{m}) ^{(2)} _{\mathbb{Q}}$ and $K_{2} ^{M} (k_{m})^{\circ} _{\mathbb{Q}},$ and the natural map from $K_{2}(k[t],(t^m))_{\mathbb{Q}}$  to $T_{m}B_{2}(k)$ obtained from the localization sequence surjects to this $K_{3} ^{\circ} (k_{m}) ^{(2)} _{\mathbb{Q}}$ as one can see by considering the reduction modulo $(t^{2m-1})$ map below.  

The reduction modulo $(t^{2m-1})$ map:
$$
(K_{2}(\pazocal{O},(t^m))/im(K_{1}(k) \cdot K_{1}(\pazocal{O},(t^m)))_{\mathbb{Q}}  \to (K_{2}(k_{2m-1},(t^m))/im(K_{1}(k) \cdot K_{1}(k_{2m-1},(t^m)))_{\mathbb{Q}},
$$
 from $T_{m}B_{2}(k) $ to $(K_{2}(k_{2m-1},(t^m))/im(K_{1}(k) \cdot K_{1}(k_{2m-1},(t^m)))_{\mathbb{Q}} \simeq K_{3} ^{\circ}(k_{m})_{\mathbb{Q}} ^{(2)} \simeq \oplus _{m<w<2m} t^w k$ is the additive dilogarithm map in this context. 
 
 If one starts with the localization sequence for the ideal $(t(1-t))$ instead of the one for $(t^m)$ above, one obtains a similar complex which computes the ordinary  weight two  motivic cohomology of $k.$ This was carried out in the fundamental work \cite{bloch-reg}.

\subsection{The additive dilogarithm as an infinitesimal dilogarithm}\label{section add dilog} 

In the first part, we describe the infinitesimal  analog of the Bloch-Wigner dilogarithm. In the second part, we explain the relation of our complex to that of Cathelineau's and that of Bloch-Esnault's. We also describe how the 4-term functional equation is related to the 5-term functional equation. 

\subsubsection{Construction of $\ell i_{m,w}$}\label{section constr of additive}

For any local $\mathbb{Q}$-algebra $A,$ we let $B_{2}(A)$ denote the $\mathbb{Q}$-space generated by $[x]$ with $x \in A^{\flat}:=\{ x| x(1-x) \in A^{\times}\}$ subject to the relations (\ref{funct-dilog}), for all $x,$ $y \in A^{\times}$ such that $(1-x)(1-y)(1-x/y) \in A^{\times}.$ We then have a complex $\Gamma_{A}(2)$ as in \textsection \ref{aomoto}.

Let $k$ be a field of characteristic 0, $k_{\infty}:=k[[t]],$ the formal power series over $k,$ and for $1 \leq m,$ $k_{m}:=k_{\infty}/(t^m). $ Recall that the Bloch-Wigner dilogarithm $D$ defines a map $B_{2} (\mathbb{C}) \to \mathbb{R}$  and is the unique measurable function, up to multiplication, with this property \cite{bloch-reg}. Its restriction of $K_{3}(\mathbb{C})^{(2)} _{\mathbb{Q}}$ is, up to a rational multiple, the Borel regulator. We have the corresponding theorem for the infinitesimal part of  $B_{2}(k_m).$  In order to describe the infinitesimal analogs of $D.$ First note that the corresponding cohomology group ${\rm H}^{1}(\Gamma _{k_m} ^{\circ}(2))$ should be $K_{3} ^{\circ}(k_{m})^{(2)}_{\mathbb{Q}}. $ This last group, by Goodwillie's theorem \cite{good}, can be expressed in terms of cyclic homology, relative to $\mathbb{Q},$ as ${\rm HC}_{2} ^{\circ}(k_m) ^{(1)}.$ There is an action,  which we denote by $\star,$ of $k^{\times}$ on $k_{m}$ such that $\lambda \in k^{\times}$ acts by sending $t$ to $\lambda \star t:=\lambda t.$ The induced action on ${\rm HC}_{2} ^{\circ}(k_m) ^{(1)}$ decomposes this group into a direct sum ${\rm HC}_{2} ^{\circ}(k_m) ^{(1)}=\oplus _{m< w <2m} k,$ with respect to the weights of the $\star$ action.  The action of $\lambda \in k^{\times}$ on $k,$  in the component of $\star$-weight   $w,$ is the one which sends $a \in k$ to $\lambda ^{w} a \in k$ (\cite{unv1}, \cite{lambda}).   This suggests that corresponding to each $\star$-weight $w$ between $m$ and $2m,$ there is a dilogarithm  $\ell i_{m,w}:B_{2}(k_{m}) \to k,$ which vanishes on the image of $B_{2}(k)$ in  $B_{2}(k_{m}) $ and induces an isomorphism between the $\star$-weight $w$ component in ${\rm HC}_{2} ^{\circ}(k_m) ^{(1)}$ and the target.   

We describe this dilogarithm as follows.  Let $\log ^{\circ}: k_{\infty} ^{\times} \to k_{\infty},$ be defined as $\log ^{\circ}(\alpha):=\log (\frac{\alpha}{\alpha(0)}).$ If $q=\sum _{0\leq i} q_i t^i \in k_{\infty}$ and $1 \leq  a $  then  $q|_{a}:=\sum _{0\leq i<a} q_i t^i 
 ,$ and $t_a(q):=q_a.$ If $u \in tk_{\infty}$ and $s(1-s) \in k^{\times},$ we let    
\begin{eqnarray}\label{formula dilog}
\ell i_{m,w}(se^{u}):=t_{w-1}(\log^{\circ}(1-se^{u|_m}) \cdot \frac{\partial u}{\partial t}\big| _{w-m})  ,
\end{eqnarray}
for $m < w< 2m.$   
 
The Bloch complex $\Gamma_{k_m}(2)$ computes the motivic cohomology of weight two over the truncated polynomial ring $k_m.$ Namely the sequence: 
\[
 \begin{CD}
 0 \to K_{3}(k_{m})^{(2)} _{\mathbb{Q}} @>>> B_{2}(k_{m})@>{\delta_m}>> \Lambda^{2} k_{m}^{\times} @>>> K_{2} ^M(k_{m})_{\mathbb{Q}} \to 0
 \end{CD}
 \]
is exact. We can state the combination of these as \cite{unv1}:

\begin{theorem*}
{\it The complex  $B_{2} ^{\circ}(k_{m})  \xrightarrow{\delta ^{\circ}} (\Lambda^2 k_{m} ^{\times})^{\circ}$ computes the infinitesimal part of the weight two motivic cohomology of $k_{m},$ the maps $\ell i_{m,w}$ satisfy the functional equation for the dilogarithm and descend to give maps from $B_{2}(k_{m})$ to $k,$ such that $\oplus _{m <w<2m} \ell i_{m,w}$ induces an isomorphism 
$$
{\rm HC}_{2} ^{\circ}(k_{m})^{(1)}  \simeq K_{3} ^{\circ}(k_m)^{(2)} _{\mathbb{Q}} \simeq ker (\delta ^{\circ}) \xrightarrow{\sim} k^{\oplus (m-1)}. 
$$}
\end{theorem*}

We sketch the main points of the proof in \cite{unv1}. 
 First we describe the map $\ell i_{m,w}$ in terms of the map $\delta.$ In order to specify the range and domain of $\delta,$ we denote the $\delta$ from $B_{2}(k_m)$ to $\Lambda ^2 k_m ^{\times}$ by the symbol $\delta_m.$ For $i<w,$ let $\ell _i: k_{w} ^{\times} \to k$ be defined by $\ell_i(\alpha):=t_i(\log ^{\circ} \alpha),$ and $\ell _i \wedge \ell _j: \Lambda ^2 k_{w} ^{\times} \to k,$ as $(\ell _i \wedge \ell _j)(a\wedge b):=\ell _i(a)\cdot  \ell _j(b)-\ell _i(b) \cdot \ell _j(a).$  
 
 The following diagram 
 $$
 \xymatrix{
 B_{2}(k_w) \ar@{->>}^{\pi_{w,m}}[d] \ar^{\delta_w}[r]& \Lambda^2 k_{w} ^{\times} \ar^{\sum_{1\leq i \leq w-m}i \cdot \ell_{w-i}\wedge \ell _i }[d]\\
  B_{2}(k_m) \ar^{\ell i_{m,w}}@{-->}[r]& k,
 }
 $$
 where $\pi_{w,m}:B_{2}(k_w) \to B_{2}(k_m)$ is the natural projection, commutes. 
This shows that $\ell i_{m,w}$ satisfies the same five term functional equation as the Bloch-Wigner dilogarithm.

Next by the stabilization theorem of Suslin ${\rm H}_{3}({\rm GL}(k_m),\mathbb{Q})={\rm H}_{3}({\rm GL}_3(k_m),\mathbb{Q}),$ and by an argument of Goncharov, we have a map ${\rm H}_{3}({\rm GL}_3(k_m),\mathbb{Q}) \to {\rm H}_{3}({\rm GL}_2(k_m),\mathbb{Q}).$ Studying the action of ${\rm GL}_{2}$ on configurations of points on $\mathbb{P}^1,$ it is easy to construct a map from ${\rm H}_{3}({\rm GL}_2(k_m),\mathbb{Q})$ to ${\rm ker}(\delta_m).$ Combining these, we obtain a map from ${\rm H}_{3}({\rm GL}(k_m),\mathbb{Q})$ to ${\rm ker}(\delta_m).$

Using Volodin's construction of $K$-theory, we can then make  Goodwillie's theorem explicit by constructing a map from ${\rm HC}_{2}^{\circ}(k_m)$ to $H_{3}({\rm GL}(k_m),\mathbb{Q}).$ Combining with the above, this gives a surjection from  ${\rm HC}_{2} ^{\circ}(k_m)$ to ${\rm ker}(\delta_m),$    Finally, by an explicit computation of   $\ell i _{m,w}$ on the image of a basis of ${\rm HC}_{2}^{\circ}(k_m)$ in ${\rm ker}(\delta_m),$ we see that $\oplus _{m<w<2m}\ell i _{m,w}$ is injective on  ${\rm HC}_{2}^{\circ}(k_m).$ This implies that the above surjection is an isomorphism and that $\oplus _{m<w<2m}\ell i _{m,w}$ is injective on it.

{\bf Example.} Using the formula (\ref{formula dilog}) above one can explicitly compute the additive dilogarithms. For example, $\ell i_{2,3}:B_{2}(k_{2}) \to k$ is given by 
$$
\ell i_{2,3}([s+at])=-\frac{a^3}{2s^2(1-s)^2}.
$$

The above theorem is an exact analog of Sydler's theorem which provides a solution to Hilbert's 3rd problem. This states that the scissors congruence class of a  three-dimensional polyhedron is completely determined by its Dehn invariant and volume. In this  context $\delta_m$ corresponds to the Dehn invariant map and $\oplus _{m<w<2m} \ell i_{m,w}$ is the sum of volumes of different $\star$-weights. When $m=2,$ this analogy gets even more precise.  In this case, there is only one dilogarithm of $\star$-weight 3, and the corresponding complex, which can be thought of as the deformation of the hyperbolic scissors complex, is analogous to the euclidean scissors congruence complex and on this complex the volume map, which is the analog of the dilogarithm, scales by the cube of the dilation factor.

\subsubsection{Comparison of $\Gamma_{k_m}(2)$ to $T_{m}\mathbb{Q}(2)(k)$ and $\beta_{2}(k)$}\label{section comparison} We first describe a  subcomplex $\overline{\Gamma}^{\circ}_{k_m}(2)$ of $\Gamma_{k_m} ^{\circ}(2).$ This is the complex whose degree 2 term is $k^{\times}\otimes (1+(t))=k^{\times} \otimes (k_{m}^{\times})^{\circ} \subseteq \Lambda ^{2}k_{m} ^{\times},$ and degree 1 term is $\delta_m  ^{-1} (k^{\times}\otimes (k_{m}^{\times})^{\circ}) \subseteq B_{2} ^{\circ}(k_m).$ We denote this last group by $\overline{B}_{2}^{\circ}(k_m).$ Then the inclusion 
 is a quasi-isomorphism from  $\overline{\Gamma}_{k_m}^{\circ}(2)$ to $\Gamma_{k_m} ^{\circ}(2)$ \cite[Proposition 6.1.2]{unv1}. In \cite[Corollary 1.4.1]{unv1}, noticing that the terms in degree 2 are the same in both of the complexes and  using the dilogarithm in degree 1 we deduce that the complexes $\overline{\Gamma} ^{\circ} _{k_m}(2)$ and $T_{m}\mathbb{Q}(2)(k)$ are isomorphic. 

Let $\tilde{\beta}_{2}(k)$ denote the $\mathbb{Q}[k^{\times}]$-module generated by $\langle a \rangle ,$ for $a \in k^{\flat},$   with the action of $\lambda \in k^{\times}$ on $\alpha$ by $\lambda \star \alpha,$ subject to the relations generated by 
\begin{align*}
   & \langle a \rangle -\langle b \rangle+a\star\langle \frac{b}{a}\rangle - (a-1)\star\langle  \frac{1-b}{1-a} \rangle=0,\\
   &(-1)\star \langle 1-a \rangle=  -\langle a \rangle \;\;\;\;\;\;\;\; {\rm and}
\;\;\;\;\;\;\;\; a \star  \langle a ^{-1} \rangle=-  \langle a \rangle.
\end{align*}

For $a \in k^{\flat},$ let  $\langle a \rangle:=a+a(1-a)t \in k_2.$ Then we have the following relations, 
 $$
 \frac{\langle b \rangle}{\langle a \rangle}=a \star  \langle \frac{b}{a} \rangle ,\;\;\; \frac{1-\langle a \rangle}{1-\langle b \rangle}=(b-1)\star \langle \frac{1-a}{1-b}\rangle, \;\; {\rm and} \;\; 1-\langle a \rangle ^{-1} =(1-a^{-1})(1-t),
 $$
 and hence $$\frac{1-\langle a \rangle ^{-1}}{1- \langle b \rangle ^{-1}}=\frac{1-a^{-1}}{1-b^{-1}}.$$  
 These imply, by the 5-term relation, that 
 $$
 [\langle a \rangle ]- [\langle b \rangle ]+a \star  [\langle \frac{b}{a} \rangle ]+(b-1)\star  [\langle \frac{1-a}{1-b} \rangle ]=0
 $$
 in $B_{2} ^{\circ}(k_2).$ Since for any $x \in k_{2} ^{\flat},$ $[1-x]=-[x]$ and $[x^{-1}]=-[x] $ in $B_{2}(k_{2}),$ we have $(-1)\star [\langle 1-a \rangle]=-[\langle a \rangle]$ and $a \star[\langle a^{-1} \rangle]=-[\langle a \rangle].$  
These relations imply that the   map that sends $\langle a \rangle$ to $[\langle a \rangle] \in B_{2} ^{\circ} (k_2)$ factors through $\tilde{\beta}_{2}(k).$ There is also  a natural surjection from  $\tilde{\beta}_{2}(k) $ to $\beta_{2}(k),$ with $\beta_{2}(k)$ defined as in \textsection \ref{4-term-section}. These maps describe the 4-term functional equation of Cathelineau as a deformation of the standard 5-term functional equation computed on special elements, where one of the terms  vanish since it has no infinitesimal part.

\section{Infinitesimal Chow dilogarithm and the infinitesimal Bloch regulator} 
In this section, we will define variants of the additive dilogarithm in order to be able to construct regulators in different settings. The first section could be thought of as removing the restriction of  considering only linear configurations when defining additive dilogarithms and is the essential step in being able to apply additive dilogarithms in an algebro geometric setting. In the second part, we describe the infinitesimal version of the Bloch regulator on curves, removing the restriction of being a curve. This is the infinitesimal version of the tame symbol construction of Deligne \cite{de}.

\subsection{Infinitesimal Chow dilogarithm} 

In this section, we construct the infinitesimal analog of the Chow dilogarithm described in \textsection \ref{section chow dilog}. The details of the construction are in \cite{unv4}. We will only consider the case of $k_2,$
 the generalization of this construction to the higher modulus case is current work in progress. The specialization of this construction to the curve $\mathbb{P}^1 _{k_2}$ and to the three linear fractional functions $1-z,$ $z$  and $1-\frac{a}{z}$ gives the additive dilogarithm $\ell i _{2,3}$ constructed in \textsection \ref{section constr of additive}, \cite[Lemma 3.5.1]{unv4}.

\subsubsection{Construction of the infinitesimal Chow dilogarithm} 

In this section, we continue to assume that $k$ is a field of characteristic 0. Let $C_{2}$ be a smooth and projective curve over $k_2.$ We do {\it not} assume that $C_{2}$ comes as a product of a curve over $k$ and $k_2.$ Let $\underline{C}_{2}$ denote the fiber of $C_{2}$ over the closed point of ${\rm Spec} (k_2).$ Given $c$ a (closed) point in $\underline{C}_{2},$ we call an element $\pi_{2,c} \in \pazocal{O}_{C_{2},c}$ in the local ring of $C_{2}$ at $c$ a {\it uniformizer}, if its reduction is a uniformizer in the local ring of $\underline{C}_2$ at $c.$ We call an element $y$ in the local ring $k(C_2)$ of $C_{2}$ at its generic point,  $\pi_{2,c}$-{\it good}, if there exists $a \in \mathbb{Z}$ such that $y=\pi_{2,c} ^a u,$ for some unit $u$ in $\pazocal{O}_{C_{2},c}.$  Note that $k(C_2)$ is an artin ring with residue field equal to the function field of $\underline{C}_{2}.$ Fix a set $\mathcal{P}_{2}:=\{\pi_{2,c}|c \in \underline{C}_{2} \}$ of uniformizers. We say that $y$ is $\mathcal{P}_{2}$-{\it good}, if it is $\pi_{2,c}$-good, for all  $c\in \underline{C}_{2}.$ 
Letting $k(C_{2},\mathcal{P}_{2})^{\times} \subseteq k(C_{2})^{\times}$ denote the group of functions which are $\mathcal{P}_{2}$-good, the infinitesimal Chow dilogarithm $\rho$ is a map
$
\rho: \Lambda^{3} k(C_{2},\mathcal{P}_{2} )^ \times \to k.
 $

In the previous section, we defined the additive dilogarithm  $\ell i_{2,3}: B_{2}(k_2) \to k,$  by 
$$
\ell i_{2,3}([s+at])=-\frac{a^3}{2s^2 (1-s)^2},
$$
and interpreted this function as the function induced by the composition
$$
(\ell _2 \wedge \ell _1)\circ\delta_{\infty}:B_{2}(k_{\infty})\to \Lambda ^2 k_{\infty} ^{\times} \to k
$$
via the canonical map $B_{2}(k_{\infty}) \to B_{2}(k_2).$ Let us denote the map $\ell _{2} \wedge \ell_1:  \Lambda ^2 k_{\infty} ^{\times} \to k$ by $\ell.$

If $\pazocal{A}/k_{\infty}$ is a smooth  algebra over $k_{\infty}$ of relative  dimension 1,  $c$ is a closed point of the spectrum of its reduction $\underline{\pazocal{A}}$ modulo $(t),$ then we call an element $\tilde{\pi}_{c}$ of the local ring $\pazocal{A}$ at $c,$ a uniformizer, if its reduction is a uniformizer in the corresponding local ring at $\underline{\pazocal{A}}.$ we have similar notions of {\it goodness} with respect to $\tilde{\pi}_{c}.$  If $\tilde{f},$ $\tilde{g}$ and $\tilde{h}$ are three functions in the local ring of $\pazocal{A}$ at the generic point of $\pazocal{A},$ which are $\tilde{\pi}_c$-good, then one can define their residue along $\tilde{\pi}_c:$   
$$
res_{\tilde{\pi}_c} (\tilde{f} \wedge \tilde{g} \wedge \tilde{h}) \in \Lambda ^2 (\pazocal{A}/(\tilde{\pi}_c))^{\times}.
 $$
As $k_{\infty}$-algebras  $\pazocal{A}/(\tilde{\pi}_c)$ is canonically isomorphic to $k'_{\infty}$ for the finite extension $k'$ of $k$ which is the residue field of $c.$  Therefore, we have a well-defined element $\ell (res_{\tilde{\pi}_c} (\tilde{f} \wedge \tilde{g} \wedge \tilde{h})) \in k'$ whose trace ${\rm Tr}_k(\ell (res_{\tilde{\pi}_c} (\tilde{f} \wedge \tilde{g} \wedge \tilde{h})))$  from $k'$ to $k$ will be essential in defining the local contribution to the Chow dilogarithm. 

In case   $C_{2}/k_{2}$ has a global lifting $\tilde{C}/k_{\infty}$ to a smooth and projective curve and $f,$ $g,$ and $h$ have global liftings $\tilde{f},$ $\tilde{g}$, and  $\tilde{h}$ to functions on $\tilde{C}$ which are good with respect to a system of uniformizers $\tilde{\mathcal{P}}:=\{\tilde{\pi}_{c}|c \in |\tilde{C}_s| \}$ on $\tilde{C}$ that  lift $\mathcal{P}_{2}$ then  
\begin{eqnarray*}
\rho(f,g,h)=\sum _{c \in |\underline{C}_2|} {\rm Tr}_k(\ell(res_{\tilde{\pi}_c} (\tilde{f} \wedge \tilde{g} \wedge \tilde{h}))).
\end{eqnarray*}
 
In general, we cannot expect such global liftings to exist. The method of defining $\rho$ is then to choose a generic lifting of the curve and arbitrary liftings of the functions and for each point of the curve to choose also local liftings of the curve together with good local liftings of the functions and then to use the residues of a 1-form which measures the defect between choosing different models. We next describe this in detail. 

The 1-form in question is defined as follows. We attach an element $\omega(\tilde{p},\hat{p},\chi) \in \Omega^{1}_{\underline{\hat{\pazocal{A}}}/k}$ to the following data:  smooth affine schemes  $\tilde{\pazocal{A}},\,  \hat{\pazocal{A}}/k_{\infty}$ of relative dimension one, an isomorphism $\chi: \tilde{\pazocal{A}}/(t^2) \xrightarrow{\sim} \hat{\pazocal{A}}/(t^2)$ and triples of functions $\tilde{p}:=(\tilde{f},\tilde{g},\tilde{h})$ in $\tilde{\pazocal{A}} ^{\times}$ and  $\hat{p}:=(\hat{f},\hat{g},\hat{h})$ in $\hat{\pazocal{A}} ^{\times},$ whose reductions modulo $t^2$ map to each other via $\chi.$ Let $\overline{\chi}:\tilde{\pazocal{A}} \xrightarrow{\sim} \hat{\pazocal{A}}$  be any lifting of $\chi,$ and $\varphi: \underline{\hat{\pazocal{A}}} \to \hat{\pazocal{A}}$ be any splitting of the canonical projection, which exist because of the smoothness assumptions.  Denote by $\overline{\varphi}$ the corresponding isomorphism $\underline{\hat{\pazocal{A}}}[[t]] \xrightarrow{\sim} \hat{\pazocal{A}}.$ 
Then we let: 
$$
\omega(\tilde{p},\hat{p},\chi):= \Omega( \overline{\varphi}^{-1}( \overline{\chi}(\tilde{p})), \overline{\varphi}^{-1}(\hat{p}) ),
$$
with $\Omega$ as below. 

Let $\tilde{q}=(\tilde{y}_{1},\tilde{y}_{2},\tilde{y_{3}})$ and $\hat{q}=(\hat{y}_{1},\hat{y}_{2},\hat{y_{3}}),$ with $\tilde{y}_i, \, \hat{y}_i \in \underline{\hat{\pazocal{A}}}[[t]] ^{\times},$ and $\hat{y}_{i}-\tilde{y}_i \in (t^2),$ for all $1 \leq i \leq 3.$ Then we can   write uniquely, $\hat{y}_i=\alpha_{0i}e^{t\alpha_{1i}+t^2\hat{\alpha}_{2i}+\cdots }$ and $\tilde{y}_i=\alpha_{0i}e^{t\alpha_{1i}+t^2\tilde{\alpha}_{2i}+\cdots },$ with $\alpha_{ji}, \hat{\alpha}_{ki},\tilde{\alpha}_{ki} \in \underline{\hat{\pazocal{A}}} ,$ for $0 \leq j \leq 1, \,2 \leq k,$ and $1\leq i \leq 3.$  We then define 
\begin{eqnarray*}
\Omega (\tilde{q},\hat{q}) := \sum _{\sigma \in S_{3}} (-1) ^\sigma \alpha_{1 \sigma(1)} ( \tilde{\alpha}_{2\sigma(3)}-\hat{\alpha}_{2\sigma(3)} ) \cdot d \log (\alpha_{0\sigma(2)}) \in \Omega^{1} _{\underline{\hat{\pazocal{A}}}/k}.
\end{eqnarray*}
The definition of $\omega(\tilde{p},\hat{p},\chi)$ is then  independent of all the choices involved.

Suppose  that  $p$ is a triple $(f,g,h)$ of functions on $C_{2}$ which are   $\mathcal{P}_{2}$-good, i.e. in $k(C_{2},\mathcal{P}_2) ^{\times}.$  In order to define $\rho(p),$ we first choose generic and local liftings of $C_{2}$ as follows.  Let $\tilde{\pazocal{A}}$ be a generic lifting of $C_{2}.$ More precisely, $\tilde{\pazocal{A}}/k_{\infty}$ is a smooth algebra together with  an isomorphism 
$
\alpha: \tilde{\pazocal{A}}/(t^2)\xrightarrow{\sim} \pazocal{O}_{C_{2},\eta}.
 $
 Let $\tilde{p}_{\eta}$ be a triple of functions in $\tilde{\pazocal{A}},$ whose reductions modulo $(t^2)$ map to the germs $p_{\eta}$ of the functions $p$ at $\eta.$ For each $c\in |C_{2}| ,$ let $\widetilde{\pazocal{B}^{\circ}_{c}}$ be a local lifting of $C_{2}$ at $c.$ In other words, $\widetilde{\pazocal{B}^{\circ}_{c}}/k_{\infty}$ is a smooth algebra together with an isomorphism $
\tilde{\gamma}_{c}: \widetilde{\pazocal{B}^{\circ}_{c}}/(t^2) \xrightarrow{\sim}  \hat{\pazocal{O}}_{C_{2},c},
$ 
from the reduction of $\widetilde{\pazocal{B}^{\circ}_{c}}$ modulo $(t^2)$ to the completion of the local ring of $C_{2}$ at $c.$ Let $\tilde{q}_{c}$ be a triple of functions on the localization of $\widetilde{\pazocal{B}^{\circ}_{c}}$ at the prime ideal $(t),$ which map to the image of $p$ via the map $\tilde{\gamma}_c$  and which are good with respect to a lift of the uniformizer on $\widetilde{\pazocal{B}^{\circ}_{c}}.$ Because of this goodness assumption on $\tilde{q}_{c},$ its residue is well-defined. We can add a term which measures the defect between the choices of the local liftings and the generic lifting and define the value of $\rho$ on $p$ as:  
\begin{eqnarray*}
\rho(p):= \sum_{c\in |C|}{\rm Tr}_k(\ell(res_c(\tilde{q}_{c}) )+ res_{c} \omega(\tilde{p}_{\eta},\tilde{q}_{c} ,  \tilde{\gamma}_{c,\eta} ^{-1} \circ \alpha_c )).
\end{eqnarray*}
It turns out that this definition is independent of all the choices involved and define a map 
from $\Lambda ^3 k(C_{2},\mathcal{P}_{2}) ^{\times}$ to $k.$ 

One can define a version of the Bloch group $B_{2}(k(C_{2},\mathcal{P}_2))$  consisting of functions which are $\mathcal{P}_{2}$-good on $C_{2}$ as in \cite[\textsection 3.3]{unv4} and define a map 
$$
\Delta:B_2 (k(C_2,\mathcal{P}_{2}))\otimes k(C_2,\mathcal{P}_{2})^{\times}  \to \Lambda^{3}k(C_2,\mathcal{P}_{2}) ^{\times},
$$ sending $[f] \otimes g$ to $(1-f)\wedge f \wedge g.$ This can then be sheafified, and using the residue map,  made into a complex which computes the motivic cohomology group $ K_{3} ^{\circ}(C_{2})^{(3)} _{\mathbb{Q}}$ as we described in \textsection \ref{section chow dilog} above, in the complex case. The infinitesimal Chow dilogarithm $\rho$ and the additive dilogarithm in the previous section can then be joined together to define a regulator from  $ K_{3} ^{\circ}(C_{2})^{(3)} _{\mathbb{Q}}$ to $k.$ We will construct and analyze this map in a future paper.

\subsubsection{Goncharov's strong reciprocity conjecture in the infinitesimal case} 

The infinitesimal Chow dilogarithm allows us to state and prove an infinitesimal version of the strong reciprocity conjecture of Goncharov \cite{gon3} with an explicit formula for the homotopy map. Let us first state the original version of the conjecture over a field which was proved recently by Rudenko \cite{rud}. 

Let $C/k$ be a smooth and projective curve over an algebraically closed field $k$ of characteristic 0. Taking the sum  of the residue maps for all  $c \in |C|,$ we obtain a commutative diagram
$$
\xymatrix{B_3(k(C)) \ar[r]  & B_2(k(C))\otimes k(C)^{\times}  \ar[r] ^-{ \Delta}\ar[d] ^{res_{|C|}}&\Lambda^{3}k(C)^{\times} \ar[d]^{res_{|C|}} \\ & B_{2}(k) \ar[r]^{\delta} & \Lambda^{2}k^{\times}.}
$$
Suslin's reciprocity theorem implies that the image of the residue map from  $\Lambda^{3}k(C)^{\times}$ to $\Lambda^{2}k^{\times}$ is in the image of $\delta.$  Goncharov's strong reciprocity conjecture states that the  residue map between the complexes above is in fact homotopic to 0 with an explicit homotopy.

In the infinitesimal setting, we start with a smooth and projective curve $C_{2}/k_{2},$ where $k$ is algebraically closed and of characteristic 0. We have the following commutative  diagram: \[
\begin{CD}
&&  B_{2}(k(C_2,\mathcal{P}_2)) \otimes k(C_2,\mathcal{P}_2)^{\times }@>{\Delta}>> \Lambda^{3} k(C_2,\mathcal{P}_{2})^{\times} \\
   && @V{\oplus res_c}VV     @V{\oplus res_c}VV\\
 && \oplus _{c \in |\underline{C}_2|} B_{2}(k_2) @>{\delta}>> \oplus _{c \in |C|} \Lambda^2 k_2 ^{\times}.
\end{CD}
\]
By \cite[Proposition 3.3.3]{unv4},  the composition $(\oplus \ell i_{2,3})\circ(\oplus res_c)$ from $B_{2}(k(C_2,\mathcal{P}_2)) \otimes k(C_2,\mathcal{P}_2)^{\times }$ to $k$ is equal to $\rho\circ\Delta.$ Then the analog of Hilbert's third problem which determines structure of $B_{2} ^{\circ}(k_2)$ in the previous section implies the following infinitesimal analog of Goncharov's strong reciprocity conjecture \cite[Theorem 3.4.4]{unv4}.

\begin{theorem*}  There is an explicit map $h: \Lambda^{3}k(C_2,\mathcal{P}_{2}) ^{\times} \to B_{2}(k_2)^{\circ}$ which makes the diagram 
$$
\xymatrix{ B_2 (k(C_2,\mathcal{P}_{2}))\otimes k(C_2,\mathcal{P}_{2})^{\times} \ar[r] ^-{ \Delta }\ar[d] ^{res_{|\underline{C}_{2}|}}&\Lambda^{3}k(C_2,\mathcal{P}_{2}) ^{\times} \ar[d]^{res_{|\underline{C}_2|}} \ar@{.>}[dl] _{h} \\  B_{2}^{\circ}(k_2) \ar[r]^{\delta^{\circ}} & (\Lambda^{2}k_2^{\times})^{\circ}}
$$
commute and has the property that $h(k_{2} ^{\times} \wedge \Lambda ^{2} k(C_2,\mathcal{P}_{2}) ^{\times})=0.$ 
\end{theorem*}

\subsubsection{An infinitesimal invariant of cycles} 

The above construction gives an infinitesimal invariant of cycles of codimension 2 in the 3 dimensional space. We briefly describe this invariant in this section and refer to \cite[\textsection 4]{unv4} for the details.  This invariant is a  generalization of the invariant defined in \cite{P1}. The approach taken in \cite{P1} for considering the infinitesimal part of the motivic cohomology of $k_2$ is to use the additive chow groups defined in \cite{addchow}, where one considers cycles on $\mathbb{A}^1 _k$ which are close to the zero cycle with  multiplicity 2 near the origin. In the approach taken here, we consider all cycles on $\mathbb{A}^{1} _k,$ but identify them if they have the same reduction modulo $(t^2).$ For the Milnor range,  the additive cycle approach is the one taken in \cite{rulling}, whereas the analog of the approach of this section is the one in \cite{ps}. 

Let $S$ denote ${\rm Spec}(k_{\infty}),$ with $s$ being the  closed and $\eta$ the generic point,   $\square_k:= \mathbb{P}^{1} _{k} \setminus \{ 1\}$ and $\square ^n _{k}$ the $n$-fold product of $\square_k$ with itself over $k, $ with the coordinate functions $y_1, \cdots, y_n.$ For  a smooth $k$-scheme $X,$ we let   $\square^n _{X} :=X \times_k \square_k ^n.$ Considering the free abelian group of  {\it admissable} cycles, the cycles which intersect each of the faces properly, of codimension $q$ on $\square^n _{X}$ for varying $n,$ one gets a complex $(\underline{z}^q(X,\cdot),\partial).$ This complex considered modulo the complex of degenerate cycles is the Bloch's cubical higher Chow complex and its cohomology groups are Bloch's higher Chow groups which compute the motivic cohomology of $X$ \cite{bloch}.

 Let $\overline{\square}_{k}:= \mathbb{P}^{1} _{k},$  $\overline{\square}_{k} ^{n},$ the $n$-fold product of  $\overline{\square}_{k} $ with itself over $k,$ and  $\overline{\square}_{S} ^{n} :=\overline{\square}_{k} ^{n} \times _k S.$ Let   $\underline{z}^q _{f} (S, \cdot) \subseteq \underline{z}^q (S, \cdot)$ be  the subgroup generated by integral, closed subschemes $Z \subseteq \square_S ^n$ which are admissible  and have {\it finite reduction}, i.e. $\overline{Z}$  intersects each $s\times \overline{F}$ properly on $\overline{\square}_S ^n,$  for  every face $F$ of $\square^n_{k}.$  Modding out by degenerate cycles, we  have the complex $z^q_{f} (S, \cdot).$

An irreducible cycle $p$ in $ \underline{z}_{f} ^2 (S,2)$ is  given by a closed point $p_{\eta} $ of  $\square ^{2} _{\eta}$ whose closure $\overline{p}$  in $\overline{\square} ^{2} _{S}$ does not meet $(\{ 0,\infty\} \times \overline{\square}_{S}) \cup ( \overline{\square}_{S} \times \{0, \infty \}).$ Let $\tilde{p}$ denote the normalisation of $\overline{p} $ and $\{s_1,\cdots, s_m \}$ the closed fiber of $\tilde{p}.$ We have  surjections $\hat{\pazocal{O}}_{\tilde{p},s_i} \to k(s_i).$ Since $k(s_i)/k$ is finite \'{e}tale there is a unique splitting $\sigma_{\tilde{p},s_i}:k(s_i)\to \hat{\pazocal{O}}_{\tilde{p},s_i}.$   We define $
\log ^{\circ} _{\tilde{p},s_i}: \hat{\pazocal{O}}_{\tilde{p},s_i} ^{\times} \to \hat{\pazocal{O}}_{\tilde{p},s_i},
$
by 
$$
\log ^{\circ} _{\tilde{p},s_i}(y)=\log(\frac{y}{\sigma_{\tilde{p},s_i}(y(s_i))}).
$$

Let 
\begin{eqnarray*}\label{defnl} 
\;\; \; l(p):=\sum _{1 \leq i \leq m}{\rm Tr}_{k}\Big(res_{\tilde{p}, s_i}\Big(\frac{1}{t^3}\big(\log^{\circ} _{\tilde{p},s_i} (y_1) \cdot d\log(y_2)-\log^{\circ} _{\tilde{p},s_i}(y_2) \cdot d \log(y_1) \big) \Big) \Big).
\end{eqnarray*}
The infinitesimal invariant 
$
\rho _f: \underline{z}_{f} ^2 (S,3) \to k
$ 
 is then defined as the composition $l  \circ \partial .$ Since $\partial ^2=0,$ it is immediate that it vanishes on boundaries. The following property is the most essential property of $\rho_f,$ which roughly states  that $\rho_f$ depends only on the reduction of $Z$ modulo $(t^2). $

 Suppose that $Z_i$ for $i=1,2$ are two irreducible cycles in $\underline{z}^{2} _{f} (S,3).$ We say that $Z_{1}$ and $Z_{2}$ are equivalent modulo $t^m$ if the following condition $(M_{m})$ holds:
 
 (i) $\overline{Z}_{i}/S$ are smooth with $(\overline{Z}_i)_{s} \cup (\cup_{j,a} |\partial _j ^{a} Z_i|) $  a strict normal crossings divisor on $\overline{Z}_i.$
  
 and more importantly 
 
 (ii) $\overline{Z}_{1}|_{t^m}=\overline{Z}_{2} |_{t^m}.$ 
 
 Then we have the following theorem \cite{unv4}: 
 
 \begin{theorem*}
{\it If $Z_{i} \in \underline{z} _{f} ^{2}(S,3),$ for $i=1,2,$  satisfy the condition $(M_{2})$ then they have the same infinitesimal regulator value: 
$$\rho_f (Z_{1})=\rho_f(Z_{2}).$$}
\end{theorem*}
Another essential property, which would justify calling $\rho_f$  a regulator, is that it vanishes on products. More precisely, if $Z \in \underline{z}_{f} ^{2}(S,3)$ and there is  $1\leq i \leq 3$  such that $y_i$ restricted to $Z|_{t^2}$ is in $k_{2} ^{\times} $ then $\rho_f(Z)=0.$ After we take the quotient with degenerate cycles, mod $(t^2)$ equivalence and boundaries, we expect $\rho_f$ to be injective, but we are very far away from  proving such a result.

\subsection{Infinitesimal Bloch regulator}

In this section, we briefly describe the infinitesimal version of the classical Bloch regulator described above in \textsection \ref{section bloch regulator}. Details of the construction will appear in \cite{unv5}. Unlike the classical case, we do not need to restrict ourselves to the case of curves. Moreover, we do not need to assume that our schemes are smooth over a truncated polynomial ring. 

We assume that $X/k$ is a scheme over a field $k$ of characteristic 0, and if $\underline{X} \hookrightarrow X$ is the scheme $X$ together with the reduced induced structure then $\underline{X}/k$ is smooth and connected, the ideal sheaf $\pazocal{I}$ of $\underline{X}$ in $X$ is square-zero, and is  locally free, as a sheaf on $\underline{X}.$ 

There is a complex which computes the infinitesimal motivic cohomology of $X$ of weight two. Namely, for a ring  $A,$ let $\Gamma_{A}(2)$ denote the complex $B_{2}(A) \to \Lambda ^{2}A^{\times},$ and if $A$ comes together with a square-zero ideal $I,$ we let   $\Gamma_{A} ^{\circ} (2)$ denote the cone of the map $\Gamma_{A}(2) \to \Gamma_{\underline{A}}(2),$ with $\underline{A}:=A/I.$ Even though $\Gamma_{A} ^{\circ} (2)$ depends on $I,$ we suppress this dependence in the notation since $I$ will be fixed in what follows. This complex is quasi-isomorphic to the complex $B_{2} ^{\circ}(A) \to (\Lambda ^{2}A^{\times})^{\circ}.$ Sheafifying this we obtain the complex of sheaves  $\Gamma_{X} ^{\circ}(2)$ on $X.$ Let  $F\Gamma_{X}^{\circ}(2) \subseteq \Gamma_{X}^{\circ}(2)$ be the subcomplex  which agrees with $\Gamma_{X}^{\circ}(2)$ in degree 1 and is the image of $\delta,$ in degree 2.  In other words, it is the subcomplex 
$$
\xymatrix{
B_{2} ^{\circ}(\pazocal{O}_{X}) \ar^{\delta}[r] \ar@{^{(}->}[d]& \delta(B_{2} ^{\circ}(\pazocal{O}_{X}) \ar@{^{(}->}[d])\\
B_{2} ^{\circ}(\pazocal{O}_{X}) \ar^{\delta}[r] & (\Lambda ^{2} \pazocal{O}_{X} ^{\times})^{\circ}
.}
$$

The analog of the Bloch regulator in this case is the following construction. We have  regulators:
$$
\rho_1: {\rm H}^2(X, \Gamma _{X} ^{\circ}(2))\to   {\rm H}^0(X, \Omega^1 _X/d \pazocal{O}_{X})^{\circ}
$$
and
$$
\rho_2: {\rm H}^2(X, F\Gamma _{X} ^{\circ}(2))={\rm ker}(\rho_1)\to   {\rm H}^1(X,D_{1}(\pazocal{O}_X)).
$$
The first regulator $\rho_1$ is defined as follows. On a  ring $A,$ with a square-zero ideal $I$ as above, we define 
$$
\log  dlog: (\Lambda ^{2} A^{\times})^{\circ} \to (\Omega ^{1} _{A}/d A)^{\circ}
$$
by sending $a \wedge b$ with $a \in (A^{\times})^{\circ}=1+I,$ and $b \in A^{\times}$ to $\log (a) dlog (b).$ The map is well-defined and vanishes on the boundaries coming from $B_{2} ^{\circ}(A).$ Therefore it can be sheafified to obtain the map $\rho_1.$ 

The more interesting part of the regulator  is $\rho_2.$ First let us give the local construcion, then we will show how to globalize this construction using a homotopy map. For the local version, we give two equivalent constructions in \cite{unv5}, one of them computational, the other one conceptual. In this survey, we only describe the computational one since it is shorter. 

Suppose that $A$ is a $k$-algebra with a square-zero ideal $I$ as above, such that $\underline{A}/k$ is smooth. The smoothness assumption implies that there is a splitting $\tau:\underline{A} \to A$ of the canonical projection. We will define a branch of the dilogarithm corresponding to this splitting. First, using the splitting $\tau$ we regard $A $  as an $\underline{A}$-algebra.  Then we  express $A$  as a quotient $B \twoheadrightarrow A$ of a  smooth $\underline{A}$-algebra $B.$ Let $\hat{B}$ denote the completion of $B$ along the kernel of this map, $\hat{\tau}$ denote  the structure map from $\underline{A}$ to $\hat{B}, $    $\hat{J}$ be  the kernel of the projection    $\hat{B} \twoheadrightarrow A,$  and $\hat{I}$ be the  inverse image of $I$ in $\hat{B.}$   Since $I^2=0,$ we have $\hat{I}^2 \subseteq \hat{J}.$ Given this presentation, the first Andr\'{e}-Quillen homology $D_{1}(A)$ of $A$ relative to $\mathbb{Q}$ is given by $D_1(A)= \ker (\hat{J}/\hat{J}^2 \xrightarrow{d} \Omega^{1} _{\hat{B}}/\hat{J}\Omega^{1}_{\hat{B}} ).$

We define a map 
$$
\ell i_{2,\tau}: \mathbb{Q}[A^{\flat}] \to D_1(A) , 
$$
by sending $[a]$ to 
$$
-\frac{1}{2}\frac{(\tilde{a}-\hat{\tau}(\underline{a}))^3}{\hat{\tau}(\underline{a})^2(\hat{\tau}(\underline{a})-1)^2} \in \ker (\hat{J}/\hat{J}^2 \xrightarrow{d} \Omega^{1} _{\hat{B}}/\hat{J}\Omega^{1}_{\hat{B}} ) ,
$$
where $\hat{\underline{a}}:=\hat{\tau}(\underline{a}),$ with $\underline{a}$ is the image of $a$ under the map $A \twoheadrightarrow \underline{A},$ and $\tilde{a}$ is any lifting  of $a \in A$ to an element in $\hat{B}.$ 
It turns out that the definition is independent of the lifting of $a$ to an element of $\hat{B}$ and that there is a natural commutative diagram corresponding to different presentations of $A$ as quotients of smooth $\underline{A}$-algebras. From the definition it is immmediate that $\ell i_{2,\tau}$ vanishes on the image of $\tau(\underline{A}^{\flat}).$ Moreover, it satisfies the 5-term relation and hence descends to a map $\ell i_{2,\tau}: B_{2}(A) \to D_{1}(A).$

In order to compare dilogarithms corresponding to  different splittings, it is necessary to restrict to the subgroup $B_{2} ^{\circ}(A)$ of $B_{2}(A).$ Again there is a more conceptual description of this homotopy map, and again we are going to take the explicit approach.

Suppose that  $A_i,$ for $i=1,\,2$ are  $k$-algebras, with  square-zero ideals $I_i,$ which are locally free $\underline{A}_i$-modules. Suppose that  $f:A_1 \to A_2$  is a $k$-algebra homomorphism and that 
$$
\tau_1: \underline{A}_1 \to A_1
$$
and 
$$
\tau_2: \underline{A}_2 \to A_2,
$$
are splittings which are {\it not} necessarily compatible with $f$. The  homotopy, that we mentioned above, in this context is a map  
$$
h_f(\tau_1, \tau_2)
: F(\Lambda^{2}A_1^{\times})^{\circ} \to D_1(A_2), 
$$
with the property that, for every $\alpha \in B_{2} ^{\circ}(A_1),$
\begin{align*}
\ell i_{2,\tau_2 } (f(\alpha))-f_{*}(\ell i_{2,\tau_1 } (\alpha))=h_f(\tau_1,\tau_2)(\delta(\alpha)).
\end{align*}
This map is a measure of  the difference between $f \circ \tau_1$ and $\tau_2  \circ\underline{f},$ where $\underline{f}: \underline{A}_1 \to \underline{A}_{2}$ is the induced map. We give the definition of $h_{f}(\tau_1,\tau_2)$ below. First, note that the map $\theta: \underline{A}_{1} \to I_{2},$ given by $\theta(a):=f(\tau_1(a))-\tau_2(\underline{f}(a))$ is an $\underline{f}$-derivation. Let $\varphi: I_{1} \to S^{2} _{\underline{A}_2} I_{2} $ be any additive  map such that for all $a \in \underline{A}_1$ and $\alpha \in I_{1}:$ 
$$
\varphi(a \alpha)=\theta(a)\otimes f(\alpha)+\underline{f}(a) \varphi(\alpha).
$$
Since $I_1$ is by assumption a locally free $\underline{A}_1$-module, such a map exists locally. Let  

$$
H_{\varphi}: (\Lambda ^{2} A_{1} ^{\times})^{\circ} \to S^{3} _{\underline{A}_2}I_{2}
$$ be the map that sends  $(1+\alpha)\wedge (1+\beta)$ with $\alpha,$ $\beta \in I_1$ to 
$
f(\alpha)\otimes \varphi (\beta)-\varphi(\alpha)\otimes f(\beta) 
$ 
and sends $(1+\alpha) \wedge \tau_{1}(a),$ with $\alpha \in I_1$ and $a \in \underline{A}_{1} ^{\flat},$ 
to 
$
-\varphi(\alpha)\otimes \frac{\theta(a)}{\underline{f}(a)}.
$ The restriction of $F(\Lambda ^2 A_{1} ^{\times})^{\circ}$ does not depend on $\varphi$ and the image lands in $D_{1}(A_{2}),$ if we use the presentation $S^{\bigcdot} _{\underline{A}_2}I_{2} \to A_2$ to compute $D_{1}(A_{2}).$ With these identifications, $h_f(\tau_1,\tau_2)$ is  the restriction of $-\frac{3}{2}H_{\varphi}$ to $F(\Lambda^2 A_{1} ^{\times} )^{\circ}.$   

This can now be used to define $\rho_2.$ Let $\{U_i \}_{i\in I}$ be an open affine cover of $X$ and   $\tau_i$ be  splittings of $\underline{U}_i \hookrightarrow U_i.$ Let $\{ a_{ij} \}_{i,j \in I} $ be local sections of $\underline{B}_2 ^{\circ}$ on $U_{ij}$ and $\{ b_{i} \}_{i \in I}$ be local sections of $F(\Lambda^{2} \pazocal{O}_{X} ^{\times})^{\circ}$  on $U_{i}$ such that $\delta (a_{ij})=b_{j}|_{U_{ij}}-b_{i}|_{U_{ij}},$ and $a_{jk}|_{U_{ijk}}-a_{ik}|_{U_{ijk}}+a_{ij}|_{U_{ijk}}=0.$ This defines an element of ${\rm H}^{2}(X,F\Gamma_{X} ^{\circ}(2)).$  

Consider the elements 
$$
\gamma_{ij}:=\ell i_{2,\tau_i}(a_{ij})+h(\tau_i,\tau_j)(b_j) \in D_{1}(\pazocal{O}_{X}(U_{ij}))),
$$ for each $i, \, j \in I.$ These define a cocycle which gives the element in ${\rm H}^{1}(X,D_{1}(\pazocal{O}_{X})),$ which is the image of $(\{ a_{ij}\}_{i, j \in I}, \{ b_{i}\}_{i \in I})$ under $\rho_2.$ This element does not depend on any of the choices made. Using Goodwillie's theorem, we also prove  in \cite{unv5} that the map $\rho_2$
 is injective. Therefore together with $\rho_1,$ they describe the motivic cohomology group ${\rm H}^{2}(X,\Gamma_{X} ^{\circ}(2))$ completely.

\section{Complements} 
In this last section, we describe some results which are very much incomplete: first the case of higher weights, then the case of characteristic $p.$ In the last section, we discuss some open problems.

\subsection{Additive polylogarithms of higher weight}\label{higher weight} In \cite{unv2}, we  constructed an analog of the single valued $n$-polylogarithms $\mathcal{L}_{n}$ of \cite{be-de}, described in \textsection \ref{section single valued l_n} above. These functions, which we denote by  $li_n$ are the the higher weight analogs of the functions defined in \textsection \ref{section add dilog}. In this higher weight case, so far we can define these functions only in modulus $(t^2),$ i.e. for $k_{2}.$ They should, of course, exist for all $k_m.$

 \begin{theorem*} \label{expressionforpoly}
 {\it For $s+at \in k_{2} ^{\flat},$ let us define
 $$
 li_{n}(s+at)=\frac{(-1)^{n}}{n!}(a/s) ^{2n-1}  \log
 (\frac{1-se^{u}}{1-s})^{(n)}(0),
 $$ 
 where the derivative is with respect to $u .$  If $A_{\bigcdot} (k_{2})$ has a comultiplication $\Delta$ such that $\Delta_{n-1,1}(\{ x\}_{n})=\{ x\}_{n-1}\otimes x$  
 for $n \geq 3$ and $\Delta_{1,1}(\{ x\}_{n})=(1-x)\wedge x$ then the above function descends to give a map $B_{n} ' (k_{2}) \to k.$} 
 \end{theorem*}
 We also proved that these infinitesimal polylogarithms satisfy the functional equations that the ordinary polylogarithms satisfy in \cite{unv2}. One should in principle be able to define such functions on all of $A_{n}(k_2),$ rather than only on the Bloch group part.

\subsection{ Partial results in characteristic $p$}\label{char p}
In the previous sections we assumed the base field $k$ to be of characteristic 0. We would expect similar constructions in characteristic $p.$ We first note that in this section we do not assume our complexes to be tensored with $\mathbb{Q}.$ Otherwise, most of the objects in question will be equal to $0.$ Therefore, for any local ring with infinite residue field $A,$ we let $B_{2}(A)_{\mathbb{Z}}$ denote the free abelian group generated by $[x],$ with $x (1-x ) \in A^{\times}$ modulo the subgroup generated by 
(\ref{funct-dilog}) with $xy(1-x)(1-y)(1-\frac{x}{y}) \in A^{\times}.$ This gives a complex $B_{2}(A)_{\mathbb{Z}} \to \Lambda ^{2} _{\mathbb{Z}}A^{\times}$ of abelian groups. We will explore this complex for $k_2,$ when $k$ is of characteristic $p.$ 

More specifically, fix $p \geq 5$ and let $\mathbb{F}$ denote an algebraic closure of the field with $p$ elements. In the following, we let  $\mathbb{F}_{m}:=\mathbb{F}[t]/(t^m).$ In particular, $\mathbb{F}_{p}=\mathbb{F}[t]/(t^p),$ and {\it not} the field with $p$ elements. 
The additive dilogarithm $\ell i_{2,3},$ with the same formula, defines a map $\ell i_{2,3}: B_{2}(\mathbb{F}_2) \to \mathbb{F}$ of $\star$-weight 3. In characteritic $p$, there is another additive dilogarithm  of $\star$-weight 1 which does {\it not} come from chracteristic 0. Recall that a finite version of the logarithm, called the $1\frac{1}{2}$-logarithm was  defined by Kontsevich \cite{kont} as: 
 
 $$  
 \pounds_{1}(s)=\sum_{1 \leq k \leq p-1}\frac{s^{k}}{k},
 $$
 for $s \in \mathbb{F}.$ This functions satisfies: 
 $\pounds_{1}(x)=-x^{p}\pounds_{1}(\frac{1}{x}),$ $
 \pounds_{1}(x)=\pounds_{1}(1-x),
$
 and 
 \begin{eqnarray*}
 \pounds_{1}(x)- \pounds_{1}(y)+x^{p} \pounds_{1}(\frac{y}{x})+(1-x)^{p} \pounds_{1}(\frac{1-y}{1-x})=0.
 \end{eqnarray*}
 Therefore $\pounds _1 ^{1/p}$ satisfies the 4-term functional equation of the entropy function. If we let 
$$
\mathfrak{Li}_{2}([s+\alpha t]):=\frac{\alpha}{s(1-s)}\sum_{1 \leq k \leq p-1}\frac{s^{k/p}}{k},
$$
then we have \cite{unv3}:

 \begin{theorem*}
{\it $\mathfrak{Li}_{2}$ descends to give a map
$$
\mathfrak{Li}_{2}:B_{2} ^{\circ}(\mathbb{F}_{2})  \to \mathbb{F}.
$$
and together with $\ell i_{2,3}$ they give the regulator from  $B_{2} ^{\circ}(\mathbb{F}_{2})  $ to $\mathbb{F} \oplus \mathbb{F}$ which gives an isomorphism  $K_{3} ^{\circ} (\mathbb{F}_{2}) \to B_{2} ^{\circ}(\mathbb{F}_{2})   \to \mathbb{F} \oplus \mathbb{F}.$}
 \end{theorem*}
 Surprisingly Kontsevich's logarithm can be obtained using $\delta$ over a truncated polynomial ring of higher modulus in an analogous manner that $\ell i_{2,3}$ can be obtained by using $\delta$ on the Bloch complex over $\mathbb{F}_3.$ However, in order to obtain $\mathfrak{Li}_{2},$ one needs to use a much higher modulus, namely one needs to lift the elements to $\mathbb{F}_p.$

Using the notation in \textsection \ref{section constr of additive}, for each $1\leq i <p,$ we have a map $\ell _{i}:\mathbb{F}^{\times} _{p} \to \mathbb{F}$ and  a commutative  diagram: 

$$
\xymatrix{
B_{2}(\mathbb{F}_p)_{\mathbb{Z}  }   \ar^{\delta_p}[r] \ar@{->>}[d] & \Lambda_{\mathbb{Z}} ^2 \mathbb{F}_{p} ^{\times} \ar^{ -\frac{1}{2} (\sum _{1 \leq i <p }a\cdot \ell _a \wedge \ell _{p-a})^{1/p}}[d] \\   B_{2}(\mathbb{F}_2)_{\mathbb{Z}  } \ar^{\mathfrak{Li}_2}[r] & \mathbb{F} ,}
$$
which expresses $\mathfrak{Li}_2$ in the manner we were looking for. 

We would  like to mention \cite{besser} for an approach to finite polylogarithms that relates them to $p$-adic polylogarithms and \cite{ev-g} for the relation of functional equations of finite polylogarithms to those of the  classical polylogarithms.

\subsection{Further problems.} As can be seen from the discussion above, the above theory is only the starting point of a general theory of infinitesimal regulators. There are many open questions, some of which will be considered in future papers. 

In the linear part of the question, the most fundamental one is that of defining the maps 
$$
vol^{\circ} _{n,w}:A_{n}(k_{m})/P_{n}(k_m) \to k
$$
of $\star$-weight $w,$ for each $(n-1)m<w<nm.$ These maps are defined above for $n=2.$ They are also defined on the subspace $B_{n} '(k_2)$ when $m=2.$ These are the analogs of the volume maps. Using these maps, one would then try to construct maps from each cohomology group of the complex $A_{\bigcdot}(k_m)$ to various $\Omega^{i} _{k}.$ One  expects, by Goodwillie's theorem and by the computation of the cyclic homology of truncated polynomial algebras that the combination of these  regulators mapping  to the direct sum of $(m-1)$-copies of $\Omega^{i} _{k}$ gives an isomorphism from the infinitesimal part of the corresponding cohomology group.  

Solving the linear part of the above problem, we expect that one could use these maps  to define regulators for smooth projective schemes $X/k_m.$ This would be the generalization of the construction of the infinitesimal Chow dilogarithm. One would have infinitesimal invariants of higher Chow groups which are expected to give all the infinitesimal invariants of the Chow groups. This last part would require significantly new ideas. 

Another main problem is to do all of the above constructions in characteristic $p.$ As we saw above in characteristic $p $ there are significantly more regulators. An essential computation in the Milnor case is done by R\"{u}lling in \cite{rulling} in the context of the additive Chow groups. In this theory, one would have to use the residue construction in the  de Rham-Witt complex rather than the ordinary de Rham complex. 

Finally, some the aspects of the construction can be done for any artin algebra over a field. This was done in the section on infinitesimal Bloch regulator in weight two. The aim would be to generalize the above to all artin algebras over a field.  For the mixed characteristic case, let us say for $W_{m}(k),$ the case of truncated Witt vectors over a perfect field $k,$ some of the regulators are of the above form. One could aim to study them using the methods above.

\end{document}